\documentclass[11pt]{article}
\usepackage{amssymb,amsmath,amsthm,amsfonts}
\usepackage{mathrsfs,amscd}
\makeatletter
\def\@fnsymbol#1{\ensuremath{\ifcase#1\or 1\or 2\fi}}
\makeatother

\def\shd{\mathcal{D}}

\def\shm{\mathcal{M}}
\def\sho{\mathcal{O}}
\def\shs{\mathcal{S}}
\newcommand{\C}{\mathbb{C}}
\newcommand{\N}{\mathbb{N}}
\newcommand{\R}{\mathbb{R}}

\newcommand{\Z}{\mathbb{Z}}
\newcommand{\K}{\mathbf{k}}
\newtheorem{theorem}{Theorem}[section]
\newtheorem*{k-s}{(K-S)-Conjecture}
\newtheorem*{k-sl}{(K-S)-Lemma}
\newtheorem{proposition}[theorem]{Proposition}
\newtheorem{lemma}[theorem]{Lemma}
\newtheorem{corollary}[theorem]{Corollary}

\theoremstyle{definition}
\newtheorem{definition}[theorem]{Definition}

\newtheorem{example}[theorem]{Example}

\begin{document}

\author{Ana Rita Martins\footnote{The research of the author was supported by
Funda\c c{\~a}o para a Ci{\^e}ncia e Tecnologia and FEDER (project
POCTI-ISFL-1-143 of Centro de Algebra da Universidade de Lisboa)
and by Funda\c c{\~a}o Calouste Gulbenkian (Programa Est{\'i}mulo
{\`a} Investiga\c c{\~a}o).}}
\title{Temperate holomorphic solutions and regularity of
holonomic D-modules on curves}
\date{\today}

\maketitle

\begin{abstract}
Let $X$ be a complex manifold. In \cite{KS2} M. Kashiwara and P.
Schapira made the conjecture that a holonomic $\shd_X$-module
$\shm$ is regular holonomic if and only if
$R\mathcal{I}{hom}_{\beta_X\shd_X}(\beta_X\shm,\sho^t_X)$ is
regular (in the sense of \cite{KS2}), the ``only if" part of this
conjecture following immediately from \cite{KS2}. Our aim is to
prove this conjecture in dimension one.
\end{abstract}

\section{Introduction}
\hspace*{\parindent}Let $X$ be a complex manifold. In \cite{KS2}
the authors defined the notions of microsupport and regularity for
ind-sheaves and applied the results to
$$Sol^t(\shm):=R\mathcal{I}{hom}_{\beta_X\shd_X}(\beta_X\shm,\sho^t_X),$$
the ind-sheaves of tempered holomorphic solutions of
$\shd_X$-modules. They proved that $$SS(Sol^t(\shm))=Char(\shm),$$
where $Char(\shm)$ denotes the characteristic variety of $\shm$,
and that $Sol^t(\shm)$ is regular if $\shm$ is regular holonomic.
In fact, M. Kashiwara and P. Schapira made the following
conjecture:
\begin{k-s}\label{k-s}
Let $\shm$ be a holonomic $\shd_X$-module. Then $\shm$ is regular
holonomic if and only if
$R\mathcal{I}{hom}_{\beta_X\shd_X}(\beta_X\shm,\sho^t_X)$ is
regular.
\end{k-s}

In this paper, we prove that, in dimension one, the regularity of
$Sol^t(\shm)$ implies the regularity of the holonomic
$\shd_X$-module $\shm$.

We start, in Section 2, with a quick review on sheaves,
ind-sheaves, microsupport and regularity for ind-sheaves and we
recall the results on the microsupport and regularity of
$Sol^t(\shm)$, proved in \cite{KS2}. We include an unpublished
Lemma of M. Kashiwara and P. Schapira that will be essential in
our proof (see the (K-S)-Lemma).

Section 3 is dedicated to the proof of the irregularity of
$Sol^t(\shm)$, when $\shm$ is an irregular holonomic $\shd$-module
on $\C$. We may reduce our proof to the case
$\shm=\shd_X^m/\shd_X^m P$, where $P$ is a matrix of differential
operators of the form $z^N\partial_z I_m+A(z),$ with $m, N\in\N$,
$I_m$ the identity matrix of order $m$ and $A$ a $m\times m$
matrix of holomorphic functions on a neighborhood of the origin.
We also show that it is enough to prove the irregularity of
$\shs^t:=\text{H}^0(Sol^t(\shm))$ and we prove that $\shs^t$ is
irregular at $(0;0)\in T^*\C$.

As an essential step we recall a classical result that gives the
characterization of the holomorphic solutions of the differential
operator $z^N\partial_z I_m+A(z)$ in some open sectors (see
Theorem \ref{T:3}). As a consequence of this result we obtain a
characterization of $\shs^t$, the ind-sheaf of temperate
holomorphic solutions of the differential operator $z^N\partial_z
I_m+A(z)$, in some open sectors (see Corollary \ref{L:2}).

Using the characterization above we find an open sector $S$, with
$0\in \overline{S}$, a small and filtrant category $\text{I}$ and
a functor $\text{I}\to D^{[a,b]}(\K_S);i\mapsto F_i$ such that
$$\sideset{``}{"}\lim_{\substack{\longrightarrow\\ i}}F_i\simeq
\shs^t|_S.$$ Moreover, we prove that for each open neighborhood
$V$ of $(0;0)$ there exist a morphism $i\to i'$ in $\text{I}$ and
a section $u$ of $\text{H}^0(\mu{hom}(F_{i},F_{i'}))$ such that
for any morphism $i'\to i''$ in $\text{I}$, denoting by $u'$ the
image of $u$ in $\text{H}^0(\mu{hom}(F_{i},F_{i''}))$, one has
supp$(u')\cap V \not\subset SS(\shs^t)$. By the (K-S)-Lemma, this
is enough to conclude the irregularity of $\shs^t$ at
$(0;0)$.\newline

\noindent{\bf {\large Acknowledgements.}} We thank P. Schapira who
suggest us to solve the (K-S)-conjecture in dimension one, using
the tools of \cite{W}, and made useful comments during the
preparation of this manuscript. We also thank M. Kashiwara and P.
Schapira for communicating us what we call the (K-S)-Lemma in
Section 2. Finally, we thank T. Monteiro Fernandes for useful
advises and constant encouragement.

\section{Notations and review}

\hspace*{\parindent}We will follow the notations in
\cite{KS2}.\newline

\noindent\textbf{Sheaves.} Let $X$ be a real $n$-dimensional
manifold. We denote by $\pi:T^*X\rightarrow X$ the cotangent
bundle to $X$. We identify $X$ with the zero section of $T^*X$ and
we denote by $\dot{T}^*X$ the set $T^*X\backslash X$.

Let $\mathbf{k}$ be a field. We denote by Mod$(\K_X)$ the abelian
category of sheaves of $\mathbf{k}$-vector spaces on $X$ and by
$D^b(\mathbf{k}_X)$ its bounded derived category.

We denote by $\R-C(\K_X)$ the abelian category of
$\R$-constructible sheaves of $\K$-vector spaces on $X$ and by
$D^b_{\R-c}(\K_X)$ the full subcategory of $D^b(\K_X)$ consisting
of objects with $\R$-constructible cohomology.

For an object $F\in D^b(\mathbf{k}_X)$, we denote by $SS(F)$ the
\emph{microsupport of $F$}, a closed $\R^+$-conic involutive
subset of $T^*X$. We refer \cite{KS3} for details.\newline

\noindent\textbf{Ind-sheaves on real manifolds.} Let $X$ be a real
analytic manifold. We denote by I$(\K_X)$ the abelian category of
ind-sheaves on $X$, that is, the category of ind-objects of the
category Mod$^c(\K_X)$ of sheaves with compact support on $X$ (see
\cite{KS1}).

Recall the natural faithful exact functor $$\iota_X:
\text{Mod}(\K_X)\to \text{I}(\K_X); F\mapsto
\sideset{``}{"}\lim_{\substack{\longrightarrow\\ U\subset\subset
X\\ U \ \text{open}}} F_U.$$ We usually don't write this functor
and identify Mod$(\K_X)$ with a full abelian subcategory of
I$(\K_X)$ and $D^b(\K_X)$ with a full subcategory of
$D^b$(I$(\K_X)$).

The category I$(\K_X)$ admits an internal hom denoted by
$\mathcal{I}{hom}$ and this functor admits a left adjoint, denoted
by $\otimes$. If $F\simeq
\underset{\underset{i}{\longrightarrow}}{\sideset{``}{"}\lim}F_i$
and $G\simeq
\underset{\underset{j}{\longrightarrow}}{\sideset{``}{"}\lim}G_j$,
then: $$\mathcal{I}{hom}(G,F)\simeq
\lim_{\substack{\longleftarrow\\
j}}\sideset{``}{"}\lim_{\substack{\longrightarrow\\
i}}\mathcal{H}{om}(G_j,F_i),$$ $$G\otimes F\simeq
\sideset{``}{"}\lim_{\substack{\longrightarrow\\
j}}\sideset{``}{"}\lim_{\substack{\longrightarrow\\ i}}(G_j\otimes
F_i).$$

The functor $\iota_X$ admits a left adjoint
$$\alpha_X:\text{I}(\K_X)\to \text{Mod}(\K_X);
F=\sideset{``}{"}\lim_{\substack{\longrightarrow\\ i}}F_i\mapsto
\lim_{\substack{\longrightarrow\\ i}}F_i.$$ This functor also
admits a left adjoint $\beta_X:\text{Mod}(\K_X)\to\text{I}(\K_X)$.
Both functors $\alpha_X$ and $\beta_X$ are exact. We refer
\cite{KS1} for the description of $\beta_X$.

Let $X$ be a real analytic manifold. We denote by $\R-C^c(\K_X)$
the full abelian subcategory of $\R-C(\K_X)$ consisting of
$\R$-constructible sheaves with compact support. We denote by
I$\R-c(\K_X)$ the category $\text{Ind}(\R-C^c(\K_X))$ and by
$D^b_{\text{I}\R-c}(\text{I}(\K_X))$ the full subcategory of
$D^b(\text{I}(\K_X))$ consisting of objects with cohomology in
I$\R-c(\K_X)$.

\begin{theorem}[\cite{KS1}]\label{T:5}
The natural functor $D^b(\mathrm{I}\R-c(\K_X))\to
D^b_{\mathrm{I}\R-c}(\mathrm{I}(\K_X))$ is an equivalence.
\end{theorem}

Recall that there is an alternative construction of
$\text{I}\R-c(\K_X)$, using Grothendieck topologies. Denote by
$\text{Op}_{X_{sa}}$ the category of open subanalytic subsets of
$X$. We may endow this category with a Grothendieck topology by
deciding that a family $\{U_i\}_i$ in $\text{Op}_{X_{sa}}$ is a
covering of $U\in \text{Op}_{X_{sa}}$ if for any compact subset
$K$ of $X$, there exists a finite subfamily which covers $U\cap
K$. One denotes by $X_{sa}$ the site defined by this topology and
by Mod$(\K_{X_{sa}})$ the category of sheaves on $X_{sa}$ (see
\cite{KS1}). We denote by $\text{Op}^c_{X_{sa}}$ the subcategory
of $\text{Op}_{X_{sa}}$ consisting of relatively compact open
subanalytic subsets of $X$.

Let $\rho:X\to X_{sa}$ be the natural morphism of sites. We have
functors
$$\text{Mod}(\K_X)\overset{\rho_*}{\underset{\rho^{-1}}{\rightleftarrows}}\text{Mod}(\K_{X_{sa}}),$$
and we still denote by $\rho_*$ the restriction of $\rho_*$ to
$\R-C(\K_X)$ and to $\R-C^c(\K_X)$.

We may extend the functor
$\rho_*:\R-C^c(\K_X)\rightarrow\text{Mod}(\K_{X_{sa}})$ to
$\text{I}\R-c(\K_X)$, by setting: $$\begin{matrix}\lambda: &
\text{I}\R-c(\K_X)&\to & \text{Mod}(\K_{X_{sa}})\\ &
\underset{\underset{i}{\longrightarrow}}{\sideset{``}{"}\lim }F_i
& \mapsto & \underset{\underset{i}{\longrightarrow}}{\lim
}\rho_*F_i.\end{matrix}$$ For $F\in\text{I}\R-c(\K_X)$, an
alternative definition of $\lambda(F)$ is given by the formula
$$\lambda(F)(U)=\text{Hom}_{\text{I}\R-c(\K_X)}(\K_U,F).$$

\begin{theorem}[\cite{KS1}]\label{T:4}
The functor $\lambda$ is an equivalence of abelian categories.
\end{theorem}
Most of the time, thanks to $\lambda$, we identify
$\text{I}\R-c(\K_X)$ with $\text{Mod}(\K_{X_{sa}})$.

We denote by $\mathcal{C}_X^\infty$ the sheaf of complex-valued
functions of class $\mathcal{C}^\infty$ and by $\shd{b}_X$ the
sheaf of Schwartz's distributions.

Let $U$ be an open subset of $X$ and let us denote
$\Gamma(U;\mathcal{C}_X^\infty)$ by $\mathcal{C}_X^\infty(U)$.

\begin{definition}
Let $f\in\mathcal{C}_X^\infty(U)$. One says $f$ has
\emph{polynomial growth} at $p\in X$ if for a local coordinate
system $(x_1,...,x_n)$ around $p$, there exist a sufficiently
small compact neighborhood $K$ of $p$ and a positive integer $N$
such that $$\sup_{x\in K\cap U} (\text{dist}(x,K\backslash
U))^N|f(x)|<\infty.$$

We say that $f$ is \emph{tempered} at $p$ if all its derivatives
have polynomial growth at $p$. We say that $f$ is tempered if it
is tempered at any point.
\end{definition}

For each open subanalytic subset $U\subset X$, we denote by
$\mathcal{C}_X^{\infty, t}(U)$ the subspace of
$\mathcal{C}_X^\infty(U)$ consisting of tempered functions and by
$\mathcal{D}{b}_X^t(U)$ the sheaf of tempered distributions on
$U$. Recall that $\mathcal{D}{b}_X^t(U)$ is defined by the exact
sequence $$0\rightarrow \Gamma_{X\backslash
U}(X;\mathcal{D}{b}_X)\rightarrow
\Gamma(X;\mathcal{D}{b}_X)\rightarrow
\mathcal{D}{b}_X^t(U)\rightarrow 0.$$

It follows from the results of Lojasiewicz \cite{L} that $U\mapsto
\mathcal{C}_X^{\infty, t}(U)$ and $U\mapsto \shd{b}_X^t(U)$ are
sheaves on the subanalytic site $X_{sa}$, hence define
ind-sheaves. We call $\mathcal{C}_X^{\infty, t}$ (resp.
$\shd{b}_X^t$) the ind-sheaf of tempered
$\mathcal{C}^\infty$-functions (resp. tempered distributions).
These ind-sheaves are well-defined in the category
Mod$(\beta_X\shd_X)$, where $\shd_X$ now denotes the sheaf of
analytic finite-order differential operators.

Let us now recall the definition of the ind-sheaf $\sho^t_X$ of
tempered holomorphic functions in a complex manifold $X$. By
definition,
$$\sho^t_X:=R\mathcal{I}{hom}_{\beta\shd_{\overline{X}}}(\beta\sho_{\overline{X}},
\shd{b}_{X_\R}^t),$$ where $\overline{X}$ denotes the complex
conjugate manifold, $X_\R$ the underlying real analytic manifold,
identified with the diagonal of $X\times\overline{X}$ and
$\shd_{\overline{X}}$ the sheaf of rings of holomorphic
differential operators of finite order over $\overline{X}$.
$\sho^t_X$ is actually an object of $D^b(\beta_X\shd_X)$ and it is
not concentrated in degree $0$ if dim $X>1$. When $X$ is a complex
analytic curve, $\sho^t_X$ is concentrated in degree $0$.
Moreover, $\sho_X$ is $\rho_*$-acyclic and $\sho^t_X$ is a
sub-ind-sheaf of $\rho_*\sho_X$.

We end this section by recalling two results of G. Morando, which
will be useful in our proof.

\begin{theorem}[\cite{G}]\label{T:1}
Let $X$ be an open subset of $\C$, $f\in \sho_\C(X)$, $f(X)\subset
Y\subset \C$. Let $U\in\mathrm{Op}^c_{X_{sa}}$ such that
$f|_{\overline{U}}$ is an injective map. Let $h\in\sho_X(f(U))$.
Then $h\circ f\in\sho_X^t(U)$ if and only if $h\in
\sho_X^t(f(U))$.
\end{theorem}

\begin{proposition}[\cite{G}]\label{P:4}
Let $p\in z^{-1}\C[z^{-1}]$ and $U\in\mathrm{Op}^c_{X_{sa}}$ with
$0\notin U$. The conditions below are equivalent.

(i) $\exp(p(z))\in \sho_X^t(U)$.

(ii) There exists $A>0$ such that $\mathrm{Re}(p(z))<A$, for all
$z\in U$.
\end{proposition}

\noindent\textbf{Microsupport and regularity for ind-sheaves.} We
refer \cite{KS2} for the equivalent definitions for the
microsupport $SS(F)$ of an object $F\in D^b(\text{I}(\K_X))$. We
shall only recall the following useful properties of this closed
conic subset of $T^*X$.

\begin{proposition}\label{P:4}
(i) For $F\in D^b(\mathrm{I}(\K_X))$, one has $SS(F)\cap
T_X^*X=\mathrm{supp}(F).$

(ii) Let $F\in D^b(\K_X)$. Then $SS(\iota_X F)=SS(F).$

(iii) Let $F\in D^b(\mathrm{I}(\K_X))$. Then
$SS(\alpha_X(F))\subset SS(F)$.

(iii) Let $F_1\to F_2\to F_3\xrightarrow{+1}$ be a distinguished
triangle in $D^b(\mathrm{I}(\K_X))$. Then $SS(F_i)\subset
SS(F_j)\cup SS(F_k)$, for $i,j,k\in\{1,2,3\}$.
\end{proposition}

Let $J$ denotes the functor
$J:D^b(\text{I}(\K_X))\to(D^b(\text{Mod}^c(\K_X)))^\wedge$ (where
$(D^b(\text{Mod}^c(\K_X)))^\wedge$ denotes the category of
functors from the $D^b(\text{Mod}^c(\K_X))^{\text{op}}$ to
\textbf{Set}) defined by:
$$J(F)(G)=\text{Hom}_{D^b(\text{I}(\K_X))}(G,F),$$ for every $F\in
D^b(\text{I}(\K_X))$ and $G\in D^b(\text{Mod}^c(\K_X))$.

\begin{definition}[\cite{KS2}]
Let $F\in D^b(\text{I}(\K_X))$, $\Lambda\subset T^*X$ be a locally
closed conic subset and $p\in T^*X$. We say that \emph{$F$ is
regular along $\Lambda$ at $p$} if there exists $F'$ isomorphic to
$F$ in a neighborhood of $\pi(p)$, an open neighborhood $U$ of $p$
with $\Lambda\cap U$ closed in $U$, a small and filtrant category
$\text{I}$ and a functor $\text{I}\to D^{[a,b]}(\K_X);i\mapsto
F_i$ such that $J(F')\simeq \underset{\underset{i\in \text{I}
}{\longrightarrow}}{\sideset{``}{"}\lim }J(F_i)$ and $SS(F_i)\cap
U\subset \Lambda$. Otherwise, we say that $F$ is irregular along
$\Lambda$ at $p$.

We say that $F$ is regular at $p$ if $F$ is regular along $SS(F)$
at $p$. If $F$ is regular  at each $p\in SS(F)$, we say that $F$
is regular.
\end{definition}

\begin{proposition}[\cite{KS2}]\label{P:1}
(i) Let $F\in D^b(\mathrm{I}(\K_X))$. Then $F$ is regular along
any locally closed set $S$ at each $p\notin SS(F)$.

(ii) Let $F_1\to F_2\to F_3\xrightarrow{+1}$ be a distinguished
triangle in $D^b(\mathrm{I}(\K_X))$. If $F_j$ and $F_k$ are
regular along $S$, so is $F_i$, for $i,j,k\in\{1,2,3\}$, $j\neq
k$.

(iii) Let $F\in D^b(\K_X)$. Then $\iota_X F$ is regular.
\end{proposition}

The next result is an unpublished Lemma of M. Kashiwara and P.
Schapira and it is very useful in the study of regularity in
$D^b(\text{I}(\K_X))$.

\begin{k-sl}
Let $\mathcal{F}\in D^b(\mathrm{I}(\K_X))$, $\Lambda\subset T^*X$
be a locally closed conic subset and  $p\in T^*X$. Assume that
there exist an open subset $S$ of $X$, with $\pi(p)\in
\overline{S}$, a small and filtrant category $I$ and a functor
$F:\mathrm{I}\to D^{[a,b]}(\K_S);i\mapsto F_i$ such that
$J(\mathcal{F}|_S)\simeq\underset{\underset{i}{\longrightarrow}}{\sideset{``}{"}\lim}J(F_i)$
and, for all open neighborhood $V$ of $p$, with $\Lambda\cap V$
closed in $V$, there exist a morphism $i\to i'$ in $\mathrm{I}$
and a section $u$ of $\mathrm{H}^0(\mu{hom}(F_{i},F_{i'}))$ such
that, for any morphism $i'\to i''$ in $\mathrm{I}$, denoting by
$u'$ the image of $u$ in $\mathrm{H}^0(\mu{hom}(F_{i},F_{i''}))$,
one has $\mathrm{supp}(u')\cap V \not\subset\Lambda$. Then
$\mathcal{F}$ is irregular along $\Lambda$ at $p$.
\end{k-sl}

\begin{proof}[\textbf{Proof.}]
We argue by contradiction. Assume that the hypothesis are
satisfied and also that $\mathcal{F}$ is regular along $\Lambda$
at $p$. Then there exist an open neighborhood $U$ of $p$ with
$\Lambda\cap U$ closed in $U$, a small and filtrant category
$\text{L}$ and a functor $G:\text{L}\to D^{[a,b]}(\K_X);l\mapsto
G_l$ such that $J(\mathcal{F}|_{\pi(U)})\simeq
\underset{\underset{l}{\longrightarrow}}{\sideset{``}{"}\lim}G_l|_{\pi(U)}$
and $SS(G_l)\cap U\subset \Lambda$, for all $l\in \text{L}$.

Since, by hypothesis, one has the following isomorphism on
$\pi(U)\cap S$:
\begin{equation}\label{I}
\sideset{``}{"}\lim_{\substack{\longrightarrow\\ i}}F_i\simeq
\underset{\underset{l}{\longrightarrow}}{\sideset{``}{"}\lim}G_l,
\end{equation}
for each $i\in \text{I}$, there exists $l(i)\in\text{L}$ and a
morphism $\rho_i:F_i\to G_{l(i)}$ on $\pi(U)\cap S$ and, for each
$l\in \text{L}$, there exist $i(l)\in \text{I}$ and a morphism
$\sigma_{l}:G_l\to F_{i(l)}$ on $\pi(U)\cap S$. Moreover, for each
$i\in \text{I}$ there exist $k\in \text{I}$ and morphisms $f:i\to
k$ and $g:i(l(i))\to k$ in $\text{I}$ such that
$F(g)\circ\sigma_{l(i)}\circ\rho_i=F(f)$.

Let $i\to i'$ be a morphism in $\text{I}$ and $u$ a section of
$\text{H}^0(\mu{hom}(F_{i},F_{i'}))$ such that for any morphism
$i'\to i''$ in $\text{I}$, denoting by $u'$ the image of $u$ in
$\text{H}^0(\mu{hom}(F_{i},F_{i''}))$, one has supp$(u')\cap U
\not\subset\Lambda$. Let $k\in \text{I}$ and $f:i'\to k$,
$g:i(l(i'))\to k$ morphisms in $\text{I}$ such that
$F(g)\circ\sigma_{l(i')}\circ\rho_{i'}=F(f)$. Then $\rho_{i'}$ and
$F(g)\circ\sigma_{l(i')}$ induce, respectively, the following
morphisms
\begin{equation}\label{E}
\text{H}^0(\mu{hom}(F_{i},F_{i'}))\to
\text{H}^0(\mu{hom}(F_{i},G_{l(i')}))\to
\text{H}^0(\mu{hom}(F_{i},F_k)),
\end{equation}
that send the section $u$ of $\text{H}^0(\mu{hom}(F_{i},F_{i'}))$
to a section $u'$ of $\text{H}^0(\mu{hom}(F_{i},$ $F_k))$.

Since supp$(\text{H}^0(\mu{hom}(F_{i},G_{l(i')})))\cap U\subset
\Lambda$, one has supp$(u')\cap U\subset\Lambda$, which is a
contradiction.
\end{proof}

\noindent\textbf{Temperate holomorphic solutions of a
$\shd$-module.} Let $X$ be a complex manifold and let $\shm$ be a
coherent $\shd_X$-module. Set
$$Sol(\shm)=R\rho_*R\mathcal{H}{om}_{\shd_X}(\shm,\sho_X),$$
$$Sol^t(\shm)=R\mathcal{I}{hom}_{\beta_X\shd_X}(\beta_X\shm,\sho^t_X).$$
It is proved in \cite{KS2} the equality:
\begin{equation}\label{E:8}
SS(Sol^t(\shm))=Char(\shm),
\end{equation}
and that the natural morphism $Sol^t(\shm)\to Sol(\shm)$ is an
isomorphism, when $\shm$ is a regular holonomic $\shd_X$-module,
which proves the ``only if" part of the (K-S)-Conjecture.

\section{Proof of (K-S)-Conjecture in dimension one}

\hspace*{\parindent}In this section, we consider $X=\C$ endowed
with the holomorphic coordinate $z$ and we shall prove that, for
every irregular holonomic $\shd_X$-module $\shm$, $Sol^t(\shm)$ is
irregular.\newline

We shall first reduce the proof to the case where
$\shm=\shd_X/\shd_X P$, for some $P\in\shd_X$.

Let $\shm$ be an irregular holonomic $\shd_X$-module. Since $\shm$
is holonomic it is locally generated by one element and we may
assume $\shm$ is of the form $\shd_X/\mathcal{I}$, for some
coherent left ideal $\mathcal{I}$ of $\shd_X$. We may also assume
that $Char(\shm)\subset T_X^*X\cup T_{\{0\}}^*X$. Moreover, we may
find $P\in\mathcal{I}$ such that the kernel of the surjective
morphism $$\shd_X/\shd_X P\to \shm\to 0,$$ is isomorphic to a
regular holonomic $\shd_X$-module $ \mathcal{N}$ (see, for
example, Chapter VI of \cite{KK}). Therefore, we have an exact
sequence $$0\to \mathcal{N}\to\shd_X/\shd_X P\to \shm\to 0,$$ and
we get the distinguished triangle $$Sol^t(\shm)\to
Sol^t(\shd_X/\shd_X P)\to Sol^t(\mathcal{N})\xrightarrow{+1}.$$
Since $Sol^t(\mathcal{N})$ is regular, by Proposition \ref{P:1},
$Sol^t(\shm)$ will be regular if and only if $Sol^t(\shd_X/\shd_X
P)$ is.

Let us now recall the following result, due to G. Morando:

\begin{theorem}[\cite{G}]
Let $\shm$ be a holonomic $\shd_X$-module. The natural morphism
$$\mathrm{H}^1(Sol^t(\shm))\rightarrow\mathrm{ H}^1(Sol(\shm)),$$
is an isomorphism.
\end{theorem}

Let $\shm=\shd_X/\shd_X P$, with $P\in\shd_X$, with an irregular
singularity at the origin. The Theorem above together with the
results in \cite{K3} entails that: $$\text{H}^1(Sol^t(\shm))\simeq
\text{H}^1(Sol(\shm))\simeq\C_{\{0\}}^m,$$ for some $m\in\N$. Then
$\text{H}^1(Sol^t(\shm))$ is regular and
$SS(\text{H}^1(Sol^t(\shm)))=T_{\{0\}}^*X$.

As in \cite{KS2}, let us set for short $$\shs^t:=
\text{H}^0(Sol^t(\shm))\simeq
\mathcal{I}{hom}_{\beta_X\shd_X}(\beta_X\shm,\sho^t_X),$$ $$\shs:=
\text{H}^0(Sol(\shm))\simeq
\rho_*\mathcal{H}{om}_{\shd_X}(\shm,\sho_X)\simeq
ker(\rho_*\sho_X\xrightarrow{P}\rho_*\sho_X).$$ Remark that, since
dim$X=1$, one has a monomorphism $\shs^t\to\shs$. Moreover, we
have the following distinguished triangle: $$\shs^t\to
Sol^t(\shm)\to \text{H}^1(Sol^t(\shm))[-1]\xrightarrow{+1}.$$
Therefore, one has $$SS(\shs^t)\subset Char(\shm)\cup
T_{\{0\}}^*X\subset T_X^*X\cup T_{\{0\}}^*X,$$ and $\shs^t$ will
be irregular if and only if $Sol^t(\shm)$ is.

The problem is then reduced to study the irregularity of $\shs^t$,
for an irregular holonomic $\shd_X$-module of the form
$\shm=\shd_X/\shd_X P$, with $P\in\shd_X$, with an irregular
singularity at the origin.

We shall prove that $\shs^t$ is not regular at $p=(0;0)$. The plan
of the proof is to find an open subset $S$ of $\C\backslash\{0\}$,
with $0\in \overline{S}$, a small and filtrant category $\text{I}$
and a functor $\text{I}\to D^{[a,b]}(\K_S);i\mapsto F_i$ such that
$$\sideset{``}{"}\lim_{\substack{\longrightarrow\\ i}}F_i\simeq
\shs^t|_S,$$ and for all open neighborhood $V$ of $p$  there exist
a morphism $i\to i'$ in $\text{I}$ and a section $u$ of
$\text{H}^0(\mu{hom}(F_{i},F_{i'}))$ such that for any morphism
$i'\to i''$ in $\text{I}$, denoting by $u'$ the image of $u$ in
$\text{H}^0(\mu{hom}(F_{i},F_{i''}))$, one has supp$(u')\cap V
\not\subset SS(\shs^t)$. This proves that we may apply the
(K-S)-Lemma to conclude that $\shs^t$ is irregular at $p$.\newline

Let $U$ be an open neighborhood of the origin in $\C$. The problem
of finding the solutions of the differential equation $Pu=0$ in
$\sho_X(U)$ is equivalent to the one of finding the solutions in
$\sho_X(U)^m$ of a system of ordinary differential equations
defined by a matrix of differential operators of the form
$$z^N\partial_z I_m+A(z),$$ where $m, N\in\N$, $I_m$ is the
identity matrix of order $m$ and $A\in
\text{M}_m(\sho_X(U))\footnote{For a ring $R$ we denote by
M$_m(R)$ the ring of $m\times m$ matrices and by $\text{GL}_m(R)$
the group of invertible $m\times m$ matrices.}$. In fact, setting
\begin{equation}\label{E:1}
P=z^N\partial_z I_m+A(z),
\end{equation}
we have $$\shs\simeq
\rho_*\mathcal{H}{om}_{\shd_X}(\shd_X^m/\shd_X^m P,\sho_X),$$
$$\shs^t\simeq
\mathcal{I}{hom}_{\beta_X\shd_X}(\beta_X(\shd_X^m/\shd_X^m
P),\sho^t_X),$$ and we may replace $\shm$ by $\shd_X^m/\shd_X^m
P$.

Let $\theta_0, \theta_1, R\in \R$, with $\theta_0<\theta_1$ and
$R>0$. We denote the open set $$\{z\in\C; \theta_0<\arg z<
\theta_1, 0<|z|<R\},$$ by $S(\theta_0,\theta_1,R)$ and call it
\emph{open sector of amplitude $\theta_1-\theta_0$ and radius
$R$}.

Let $l\in\N$. We consider $z^{1/l}$ as a holomorphic function on
subsets of open sectors of amplitude smaller than $2\pi$, by
choosing the branch of $z^{1/l}$ which has positive real values on
$\arg z=0$.

The next Theorem will be fundamental in our proof. It gives a
characterization of the holomorphic solutions of the differential
system $Pu=0$, where $P$ is the operator (\ref{E:1}).

\begin{theorem}[see \cite{W}]\label{T:3}
Let $P$ be the matrix of differential operators (\ref{E:1}). There
exist $l\in\N$ and a diagonal matrix $\Lambda(z)\in
\mathrm{M}_m(z^{-1/l}\C[z^{-1/l}])$ such that, for any real number
$\theta$, there exist $R>0$, $\theta_1>\theta>\theta_0$ and
$F_{\theta}\in \mathrm{GL}_m(\sho_X(S(\theta_0,\theta_1,R))\cap
C^0(\overline{S(\theta_0,\theta_1,R)}\backslash\{0\}))$, such that
the $m$-columns of the matrix $F_\theta(z)\exp(-\Lambda(z))$ are
$\C$-linearly independent holomorphic solutions of the system
$Pu=0$. Moreover, for each $\theta$ there exist constants $C,M>0$
so that $F_\theta$ has the estimate
\begin{equation}\label{E:2}
C^{-1}|z|^M<|F_\theta(z)|<C|z|^{-M}, \ \text{for any $z\in
S(\theta_0,\theta_1,R)$}.
\end{equation}
\end{theorem}

If there is no risk of confusion we shall write $F(z)$ instead of
$F_\theta$.

\begin{definition}
We call the matrix $F(z)\text{exp}(-\Lambda(z))$, given in Theorem
\ref{T:3}, a \emph{fundamental solution of $P$ on
$S(\theta_0,\theta_1,R)$}.
\end{definition}

Let $l\in\N$ and $\Lambda(z)$ be the diagonal matrix given in
Theorem \ref{T:3}. For each $1\leq j\leq m$, let
$\Lambda_j(z)=\sum_{k=1}^{n_j}a^j_kz^{-k/l}$ be the $(j,j)$ entry
of $\Lambda(z)$, with $n_j\in\N$, $a^j_1,...,a^j_{n_j}\in \C$.

\begin{corollary}\label{L:2}
Let $V\in\mathrm{Op}^c_{X_{sa}}$ and let us suppose $P$ has a
fundamental solution $F(z)\exp(-\Lambda(z))$ on $V$. Then,
$\Gamma(V;\shs^t)\simeq \C^{n(V)}$, where $n(V)$ is the
cardinality of the set: $$J(V):= \{j\in\{1,...,m\};
\exp(-\Lambda_j(z))|_V \in \sho^t_{X}(V)\}.$$
\end{corollary}

\begin{proof}[\textbf{Proof.}]
By hypothesis, $\Gamma(V;\shs)$ is the $m$-dimensional $\C$-vector
space generated by the $m$-columns of the matrix
$F(z)\exp(-\Lambda(z))$. For each $j=1,...,m$, let us denote by
$e_j$ the $j$-th column of the matrix $F(z)\exp(-\Lambda(z))$. We
have:
$$\Gamma(V;\shs^t)=\Gamma(V;\shs)\cap\Gamma(V;\sho_{X}^{t,m})
=\{u\in \Gamma(V;\sho_{X}^{t,m});$$ $$u=F(z)\exp(-\Lambda(z))C, \
\text{for some scalar column matrix $C$}\}.$$

Let $k$ be the dimension of the $\C$-vector space
$\Gamma(V;\shs^t)$ and let us prove that $k=n(V)$.

For each $j=1,...,m$, we have $e_j=F(z)\exp(-\Lambda(z))C_j$,
where $C_j$ is the scalar column matrix with the $j$-th entry
equal to $1$ and all the other entries equal to zero. Since $F(z)$
is a matrix of tempered holomorphic functions on $V$, we get
$e_j\in\Gamma(V;\shs^t)$, for each $j\in J(V)$, and the
$\C$-vector space generated by the family $\{e_j\}_{j\in J(V)}$ is
a vector subspace of $\Gamma(V;\shs^t)$. It follows that $n(V)\leq
k$. On the other hand, we may find a subset $K$ of $\{1,...,m\}$,
with cardinality $k$, such that $e_j\in\Gamma(V;\shs^t)$, for all
$j\in K$. Since $F^{-1}(z)$ is a matrix of tempered holomorphic
functions on $V$, we get $\exp(-\Lambda_{j}(z))\in
\sho_{X}^{t}(V)$, for each $j\in K$. This entails that $n(V)\geq k$
and completes our proof.
\end{proof}

\begin{lemma}\label{L:9}
Let $S$ be an open sector of amplitude smaller than $2\pi$, $p\in
z^{-1}\C[z^{-1}]$, $l\in\N$ and $V\in\mathrm{Op}^c_{X_{sa}}$, with
$V\subset S$. Then $\exp(p(z^{1/l}))\in\sho^t_{X}(V)$ if and only
if there exists $A>0$ such that $\mathrm{Re}(p(z^{1/l}))<A$, for
all $x\in V$.
\end{lemma}

\begin{proof}[\textbf{Proof.}]
Let $\theta_0,\theta_1, R\in\R$ such that $0<\theta_1-\theta_0<R$
and $S=S(\theta_1,\theta_0,R)$, and let us denote by $U$ the open
sector $S(\frac{\theta_0}{l},\frac{\theta_1}{l},R^{1/l})$. Let
$f:X\to X$ be the holomorphic function defined by $f(z)=z^l$.
Since $\theta_1-\theta_0<2\pi$, we may easily check that
$f|_{\overline{U}}$ is an injective map. Moreover, $f(U)=S$ and
$f|_U:U\to S$ is bijective. Set $V'=f^{-1}(V)\cap U$ and let $h$
denotes the holomorphic function defined for each $z\in S$ by
$h(z)=\exp(p(z^{1/l}))$. By Theorem \ref{T:1}, we have $h\circ
f\in \sho^t_{X}(V')$ if and only if $h\in \sho^t_{X}(V)$. On the
other hand, one has $p|_{V'}=h\circ f|_{V'}$ and, by Proposition
\ref{P:4}, $h\circ f\in \sho^t_{X}(V')$ if and only if there
exists $A>0$ such that $\text{Re}(p(z))<A$, for all $z\in V'$.
Combining these two facts, we conclude that
$\exp(p(z^{1/l}))\in\sho^t_{X}(V)$ if and only if there exists
$A>0$ such that $\text{Re}(p(z^{1/l}))<A$, for all $z\in V$, as
desired.
\end{proof}

\begin{proposition}\label{P:2}
With the notation above, there exist an open sector $S$, with
amplitude smaller than $2\pi$ and radius $R>0$, and a non-empty
subset $\mathrm{I}$ of $\{1,...,m\}$ such that, for each $j\in
\mathrm{I}$ and each open subanalytic subset $V\subset S$, the
conditions below are equivalent:

(i) there exists $A>0$ such that $\mathrm{Re}(-\Lambda_j(z))<A$
for all $z\in V$,

(ii) there exists $0<\delta<R$ such that $V\subset\{z\in
S;|z|>\delta\}$.

Moreover, for each $j\in \{1,...,m\}\backslash \mathrm{I}$, there
exists $A>0$ such that, for every $z\in S$,
$\mathrm{Re}(-\Lambda_j(z))<A$.
\end{proposition}

We shall need the following Lemma:

\begin{lemma}\label{L:3}
Let $m,l\in\N$, $\phi_1,...,\phi_m\in [0,2\pi[$ and
$n_1,...,n_m\in\N$. For each $j=1,...,m$, $0<C<1$ and
$\theta\in\R$, let us consider the following two conditions:

$(i)_{j,C,\theta}$ $1\geq\cos(\phi_j-n_j/l\theta)\geq C$,

$(ii)_{j,C,\theta}$ $-C\geq\cos(\phi_j-n_j/l\theta)\geq -1$.

Then we may find $\theta_0,\theta_1\in\R$, with
$0<\theta_1-\theta_0< 2\pi$, and positive real numbers
$C_1,...,C_m$ such that, for each $j=1,...,m$, one of the two
conditions $(i)_{j,C_j,\theta}$ or $(ii)_{j,C_j,\theta}$ holds,
for every $\theta_0\leq\theta\leq\theta_1$. Moreover, condition
$(i)_{j,C_j,\theta}$ is satisfied, for every
$\theta_0\leq\theta\leq\theta_1$, for some $j=1,...,m$.
\end{lemma}

\begin{proof}[\textbf{Proof.}]
Let us prove the result by induction on $m$. If $m=1$, set
$C_1=1/2$, $\theta_0=0$ and $\theta'_1=\frac{l\pi}{4n_1}$ if
$\phi_1=0$, and
$\theta_0=\max\{0,\frac{l}{n_1}(\phi_1-\frac{\pi}{4})\}$,
$\theta'_1=\frac{l\phi_1}{n_1}$ if $\phi_1\neq 0$. Then, for every
$\theta_0\leq\theta\leq\theta'_1$, we have
$1\geq\cos(\phi_1-n_1/l\theta)\geq C_1$. Since $0\leq\theta_0<
2l\pi$, we may find $0\leq k < l$ such that
$2k\pi\leq\theta_0<2(k+1)\pi$ and $\theta_1\in\R$ such that
$\theta_0<\theta_1<\theta'_1$ and $\theta_1<2(k+1)\pi$. Therefore,
condition $(i)_{1,C_1,\theta}$ is satisfied, for every
$\theta_0\leq\theta\leq\theta_1$, where
$0<\theta_1-\theta_0<2\pi$.

Let us now assume that the result is true for some $m\geq 1$ and
let us consider $\phi_1,...,\phi_{m+1}\in [0,2\pi[$ and
$n_1,...,n_{n+1}\in\N$. By hypothesis, there exist
$\theta'_0,\theta'_1\in\R$, with $0<\theta'_1-\theta'_0< 2\pi$,
and positive real numbers $C_1,...,C_m$ such that, for each
$j=1,...,m$, one of the two conditions $(i)_{j,C_j,\theta}$ or
$(ii)_{j,C_j,\theta}$ holds, for every
$\theta'_0\leq\theta\leq\theta'_1$. Moreover, condition
$(i)_{j,C_j,\theta}$ is satisfied, for every
$\theta'_0\leq\theta\leq\theta'_1$, for some $j=1,...,m$. Let us
prove that there exist $\theta_0,\theta_1\in\R$, with
$\theta'_0<\theta_0<\theta_1<\theta'_1$, and $0<C_{m+1}<1$ such
that one of the two conditions $(i)_{m+1, C_{m+1}, \theta}$ or
$(ii)_{m+1, C_{m+1}, \theta}$ holds, for every
$\theta_0\leq\theta\leq\theta_1$.

For each $j\in\Z$, set: $$I_j=]2j\pi,\pi/2+2j\pi[, \
J_j=]3\pi/2+2j\pi,2\pi(1+j)[,$$ $$K_j=]\pi/2+2j\pi,\pi+2j\pi[, \
L_j=]\pi+2j\pi,3\pi/2+2j\pi[.$$

If
$[\phi_{m+1}-n_{m+1}/l\theta'_1,\phi_{m+1}-n_{m+1}/i\theta'_0]\subset
A_j$, for some $A\in\{I, J, K, L\}$ and $j\in\Z$, set:
$$C_{m+1}=\left\{
\begin{array}{l}\ \ \ \cos(\phi_{m+1}-n_{m+1}/l\theta'_0),\ \text{if $A=I,$}\\
\ \ \ \cos(\phi_{m+1}-n_{m+1}/l\theta'_1),\ \text{if $A=J$,}\\
-\cos(\phi_{m+1}-n_{m+1}/l\theta'_1),\ \text{if $A=K$.}\\
-\cos(\phi_{m+1}-n_{m+1}/l\theta'_0),\ \text{if $A=L,$}
\end{array}
\right.$$ Then  condition $(i)_{m+1, C_{m+1}, \theta}$ holds, for
all $\theta'_0\leq\theta\leq\theta'_1$, if $A\in\{I, J\}$ and
condition $(ii)_{m+1, C_{m+1}, \theta}$ holds, for all
$\theta'_0\leq\theta\leq\theta'_1$, if $A\in\{K, L\}$.

If
$[\phi_{m+1}-n_{m+1}/l\theta'_1,\phi_{m+1}-n_{m+1}/l\theta'_0]\not\subset
A_j$, for every $A\in\{I, J, K, L\}$ and $j\in\Z$, we may choose
$\theta_0,\theta_1\in\R$ such that
$\theta'_0<\theta_0<\theta_1<\theta'_1$ and
$[\phi_{m+1}-n_{m+1}/l\theta_1,\phi_{m+1}-n_{m+1}/l\theta_0]\subset
A_j$, for some $A\in\{I, J, K, L\}$ and $j\in\Z$. Then the result
follows from the previous cases.
\end{proof}

\begin{proof}[\textbf{Proof of Proposition \ref{P:2}}]
For each $j=1,...,m$, if $z=\rho \exp(i \theta)$, one has:
$$\text{Re}(-\Lambda_j(z))=\sum_{k=1}^{n_j}\alpha^j_k\rho^{-k/l}\cos(\phi^j_k-k/l\theta),$$
where $\alpha^j_k=|a_k^j|$ and $\phi^j_k=\arg(-a_k^j)$, for every
$k=1,...,n_j$.

For each $j=1,...,m$, $1>C>0$ and $\theta\in\R$, let us consider
the following two conditions:

$(i)_{j,C,\theta}$ $1\geq\cos(\phi^j_{n_j}-n_j/l\theta)\geq C$,

$(ii)_{j,C,\theta}$ $-C\geq\cos(\phi^j_{n_j}-n_j/l\theta)\geq -1$.

By Lemma \ref{L:3}, we may find $\theta_0,\theta_1\in\R$, with
$0<\theta_1-\theta_0<2\pi$, and positive real numbers
$C_1,...,C_m$ such that, for each $j=1,...,m$, one of the two
conditions $(i)_{j,C_j,\theta}$ or $(ii)_{j,C_j,\theta}$ holds,
for every $\theta_0\leq\theta\leq\theta_1$. Moreover, condition
$(i)_{j,C_j,\theta}$ is satisfied, for every
$\theta_0\leq\theta\leq\theta_1$, for some $j=1,...,m$.

Let us set: $$\text{J}:=\{ j\in\{1,...,m\}; \ \text{condition
$(ii)_{j, C_j, \theta}$ is satisfied, for every
$\theta_0\leq\theta\leq\theta_1$}\}.$$ For each $j\in \text{J}$,
$\theta_0\leq\theta\leq\theta_1$ and $\rho>0$, one has:
$$\text{Re}(-\Lambda_j(\rho\exp(i\theta)))=$$
$$=\rho^{-n_j/l}[\sum_{k=1}^{n_j-1}\alpha^j_k\rho^{(n_j-k)/l}\cos
(\phi^j_k-k/l\theta)+\alpha^j_{n_j}\cos(\phi^j_{n_j}-n_j/l\theta)]\leq$$
$$\leq\rho^{-n_j/l}[\sum_{k=1}^{n_j-1}\alpha^j_k\rho^{(n_j-k)/l}-\alpha^j_{n_j}C_j],$$
and $$\lim_{\substack{\rho\to 0^+}}
\rho^{-n_j/l}[\sum_{k=1}^{n_j-1}\alpha^j_k\rho^{(n_j-k)/l}-\alpha^j_{n_j}C_j]=-\infty.$$
It follows that there exists $0<R_j$ such that
$\text{Re}(-\Lambda_j(\rho\exp(i\theta)))<0$, for every
$0<\rho<R_j$ and $\theta_0\leq\theta\leq\theta_1$. Therefore,
setting $R=\min\{R_j; j\in \text{J}\},$ one gets that
$\text{Re}(-\Lambda_j(z))<A$, for every $A>0$, $z\in
S(\theta_0,\theta_1,R)$ and $j\in \text{J}$.

Let us now set $$\text{I}:=\{j\in\{1,...,m\};\text{condition
$(i)_{j, C_j, \theta}$ is satisfied, for every
$\theta_0\leq\theta\leq\theta_1$}\}.$$ Let $j\in \text{I}$ and $V$
be an open subanalytic subset of the sector
$S(\theta_0,\theta_1,R)$. Suppose that there exists $A>0$ such
that $\text{Re}(-\Lambda_j(z))<A$ for every $z\in V$ and that, for
each $0<\delta<R$, there exists $z_\delta\in V$ with
$|z_\delta|\leq\delta$. For each $0<\delta<R$, let us denote:
$\rho_\delta=|z_\delta|$ and $\theta_\delta=\arg(z_\delta)$. The
sequence $\{\rho_{\delta}\}_\delta$ converges to $0$ and since
$\{\theta_{\delta}\}_\delta$ is a bounded sequence it admits a
convergent subsequence. Replacing these two sequences by
convenient subsequences, we may assume that
$\{\theta_{\delta}\}_\delta$ converges to some $\theta_2\in
[\theta_0,\theta_1]$. Then: $$\lim_{\substack{\delta\to
0^+}}\text{Re}(-\Lambda_j(\rho_{\delta}\exp(i\theta_{\delta})))=$$
$$=\lim_{\substack{\delta\to 0^+}}
\rho_{\delta}^{-n_j/l}[\sum_{k=1}^{n_j-1}\alpha^j_k\rho_{\delta}^{(n_j-k)/l}\cos
(\phi^j_k-k/l\theta_{\delta})+\alpha^j_{n_j}\cos(\phi^j_{n_j}-n_j/l\theta_{\delta})]\geq$$
$$\geq\lim_{\substack{\delta\to 0^+}}
\rho_{\delta}^{-n_j/l}[-\sum_{k=1}^{n_j-1}\alpha^j_k\rho_{\delta}^{(n_j-k)/l}+\alpha^j_{n_j}C_j]=+\infty,$$
which contradicts the fact that
$\text{Re}(-\Lambda_j(\rho_{\delta}\exp(i\theta_{\delta})))<A$,
for every $0<\delta<R$. Conversely, if $V$ is an open subanalytic
subset of $\{z\in S(\theta_0,\theta_1,R);|z|>\delta\}$, for some
$0<\delta<R$, then $V$ is contained on the compact set $\{z\in\C;
\theta_0\leq \arg z \leq\theta_1, \delta\leq|z|\leq R\}$, and
$\text{Re}(-\Lambda_j(z))$ is obviously bounded on $V$. We
conclude that $\text{I}$ is the desired subset of $\{1,...,m\}$,
with $\{1,...,m\}\backslash \text{I}=\text{J}$.
\end{proof}

Let $S(\theta'_0,\theta'_1,R')$ and $\text{I}$ be, respectively,
the open sector and the subset of $\{1,...,m\}$ given by
Proposition \ref{P:2} and let us choose $\theta_0,\theta_1,
R\in\R$, with $\theta'_0<\theta_0<\theta_1<\theta'_1$ and $R>0$,
such that the matrix of differential operators $P$ admits a
fundamental solution $F(z)\exp(-\Lambda(z))$ on the open sector
$S(\theta_0,\theta_1,R)$. Let us denote
$S=S(\theta_0,\theta_1,R)$, $S_\delta=\{z\in S; |z|>\delta\}$, for
each $0<\delta<R$, and let $n$ be the cardinality of the set
$\text{I}$. Remark that, $\text{I}\neq\emptyset$ and so, $n>0$.

\begin{proposition}\label{P:3}
One has the following isomorphism on $S$,
\begin{equation}\label{E:3}
\sideset{``}{"}\lim_{\substack{\longrightarrow\\
R>\delta>0}}(\C^{n}_{S_\delta}\oplus\C_{S}^{m-n})\xrightarrow{\sim}
\mathcal{I}{hom}_{\beta_X\shd_X}(\beta_X\shm,\sho^t_X)|_S.
\end{equation}
\end{proposition}

\begin{proof}[\textbf{Proof.}]
By Lemma \ref{L:9}, for each $j=1,...,m$ and
$V\in\text{Op}^c_{X_{sa}}$, with $V\subset S$,
$\exp(-\Lambda_j(z))\in\sho^t_X(V)$ if and only if there exists
$A>0$ such that $\text{Re}(-\Lambda_j(z))<A$ for each $z\in V$.
Let $V\in\text{Op}^c_{X_{sa}}$, with $V\subset S$. Thus, by
Proposition \ref{P:2}, for all $j\in \{1,...,m\}\backslash
\text{I}$, one has $\exp(-\Lambda_j(z))\in\sho^t_X(V)$ and, for
$j\in \text{I}$ one has $\exp(-\Lambda_j(z))\in\sho^t_X(V)$ if and
only if $V\subset S_\delta$, for some $0<\delta<R$. Therefore, by
Corollary \ref{L:2}, either $V\subset S_\delta$, for some
$0<\delta<R'$, and in this case $\Gamma(V;\shs^t)\simeq \C^m$, or
else $\Gamma(V;\shs^t)\simeq \C^{m-n}$. By Theorem \ref{T:4}, we
get the desired isomorphism.
\end{proof}

For each $0<\delta<R$, set
$F_\delta=\C^{n}_{S_\delta}\oplus\C_{S}^{m-n}$ and let us prove
that, for every open neighborhood $U$ of $p=(0;0)$,
$SS(F_\delta)\cap U\not\subset SS(\shs^t)$.

Let us set $z=x+iy$. For each $0<\delta<R$, we have,
$$SS(F_\delta)=SS(\C_{S_\delta})\cup T_S^*S= \{(x,y;0)\in T^*S
;x^2+y^2\geq \delta^2\}\cup$$ $$\cup\{(x,y;\lambda x,\lambda y
)\in T^*S; \lambda<0, x^2+y^2=\delta^2\}\cup T_S^*S.$$ Hence, for
each open neighborhood $U$ of $p$, we may find $0<\delta<R$ and
$(x,y)\in S$ such that $x^2+y^2=\delta^2$ and $(x,y;-x,-y)\in
(SS(F_\delta)\cap U)\backslash SS(\shs^t)$.

To finish the proof that $\shs^t$ is irregular at $p$, let us
recall that, for each $0<\delta<R$, one has the following natural
morphisms: $$\text{Hom}(F_{\delta},F_{\delta})\simeq
\text{H}^0(X;R\mathcal{H}{om}(F_{\delta},F_{\delta}))\simeq$$
$$\simeq \text{H}^0(T^*X;\mu{hom}(F_{\delta},F_{\delta}))\to
\Gamma(T^*X;\text{H}^0(\mu{hom}(F_{\delta},F_{\delta}))),$$ and we
shall denote by $u_\delta$ the image of
$\text{id}_{F_{\delta}}\in\text{Hom}(F_{\delta},F_{\delta})$ in
$\Gamma(T^*X;$ $\text{H}^0(\mu{hom}(F_{\delta},F_{\delta})))$.

For each $0<\varepsilon<\delta<R$, $S_\delta$ is an open subset of
$S_\varepsilon$ and we have an exact sequence
\begin{equation}\label{E:6}
0\to F_\delta \to F_\varepsilon \to F_{\delta,\varepsilon}\to 0,
\end{equation}
where $F_{\delta,\varepsilon}$ denotes the sheaf
$\C^{m-n}_{S_\varepsilon\backslash S_\delta}\oplus\C_{S}^{n}$.

Applying the functor $\mu{hom}(F_{\delta},\cdot)$ to the exact
sequence (\ref{E:6}), we obtain the distinguished triangle:
\begin{equation}\label{E:7}
\mu{hom}(F_{\delta},F_\delta)\to\mu{hom}(F_{\delta},F_\varepsilon)
\to\mu{hom}(F_{\delta},F_{\delta,\varepsilon})\xrightarrow{+1}.
\end{equation}

Since
$\mu{hom}(F_{\delta},F_\delta),\mu{hom}(F_{\delta},F_\varepsilon),
\mu{hom}(F_{\delta},F_{\delta,\varepsilon})\in D^{\geq
0}(\K_{T^*X})$, it follows from (\ref{E:7}) that we have the exact
sequence $$0\to \text{H}^0(\mu{hom}(F_{\delta},F_{\delta}))\to
\text{H}^0(\mu{hom}(F_{\delta},F_{\varepsilon})).$$ Hence,
supp$(u_\delta)=$supp$(u')$, where $u'$ is the image of
$u_{\delta}$ in
$\text{H}^0(\mu{hom}(F_{\delta},F_{\varepsilon}))$. Moreover, by
Corollary 6.1.3 of \cite{KS3}, one has
supp$(u_{\delta})=SS(F_{\delta})$. This proves that we may apply
the (K-S)-Lemma to conclude that $\shs^t$ is irregular at $p$.

We finish with an example.

\begin{example}
Let us consider the $\shd_X$-module
$$\shm=\shd_X/\shd_X(z^2\partial_z+1).$$ In this case, $\exp(1/z)$
is a fundamental solution of the differential operator
$z^2\partial_z+1$ in $X\backslash\{0\}$. Arguing as in the proof
of Proposition \ref{P:2}, we find $R>0$ with the following
property: given an open subanalytic subset $V$ of the sector
$S=S(0,\pi/4,R)$, then there exists $A>0$ such that
$\text{Re}(-1/z)<A$, for every $z\in V$, if and only if there
exists $0<\delta< R$ such that $V\subset\{z\in S;|z|>\delta\}$.
Moreover, by Proposition \ref{P:3}, one has the isomorphism below:
$$\sideset{``}{"}\lim_{\substack{\longrightarrow\\
R>\delta>0}}\C_{S_\delta}\xrightarrow{\sim}
\mathcal{I}{hom}_{\beta_X\shd_X}(\beta_X\shm,\sho^t_X)|_S.$$

In \cite{KS2}, M. Kashiwara and P. Schapira proved the following
isomorphism on $X$,
$$\sideset{``}{"}\lim_{\substack{\longrightarrow\\
\varepsilon>0}}\C_{U_\varepsilon}\xrightarrow{\sim}
\mathcal{I}{hom}_{\beta_X\shd_X}(\beta_X\shm,\sho^t_X),$$ where
$U_\varepsilon=X\backslash
\overline{B_\varepsilon(0,\varepsilon)},$ and
$B_\varepsilon(0,\varepsilon)$ denotes the open ball with center
at $(0,\varepsilon)$ and radius $\varepsilon$, for every
$\varepsilon>0$.

Let us check that
$$\sideset{``}{"}\lim_{\substack{\longrightarrow\\
R>\delta>0}}\C_{S_\delta}\simeq
\sideset{``}{"}\lim_{\substack{\longrightarrow\\
\varepsilon>0}}\C_{U_\varepsilon\cap S}.$$

It is enough to prove that for every $0<\delta<R$, there exists
$\varepsilon>0$ such that $S_\delta\subset U_\varepsilon\cap S$
and that for every $\varepsilon>0$ there exists $0<\delta<R$ such
that $U_\varepsilon\cap S\subset S_\delta$.

For each $0<\delta<R$, we have $S_\delta\subset U_{\delta/4}\cap
S$. In fact, given $x+iy\in S_\delta$ we have $x^2+y^2> \delta^2$
and $x>y>0$. It follows that $2x^2>\delta^2$ and hence, $x^2+y^2>
x^2>2x(\delta/4)$, i.e., $(x-\delta/4)^2+y^2>(\delta/4)^2$.
Conversely, given $\varepsilon>0$ and $x+iy\in U_\varepsilon\cap
S$, we have $x^2+y^2>2x\varepsilon$ and $x>y>0$. This gives
$x>\varepsilon$. Hence, $x^2+y^2>\varepsilon^2$ and  $x+iy\in
S_\varepsilon$.
\end{example}

{\small Ana Rita Martins

Centro de {\'A}lgebra da Universidade de Lisboa, Complexo 2,

Avenida Prof. Gama Pinto,

1699 Lisboa Portugal

arita@mat.fc.ul.pt}

\end{document}